\documentclass[11pt]{amsart}
\usepackage{amsmath, amsthm, amsfonts, amssymb}
\usepackage{mathrsfs}
\usepackage{bbm}
\usepackage{ifthen}
\numberwithin{equation}{section}
\usepackage{a4}

\newtheorem{thm}{Theorem}[section]
\newtheorem{cor}[thm]{Corollary}
\newtheorem{prop}[thm]{Proposition}
\newtheorem{lem}[thm]{Lemma}
\theoremstyle{remark}
\newtheorem{rem}[thm]{Remark}
\newtheorem{exmp}[thm]{Example}
\newtheorem{defn}[thm]{Definition}

\newcommand{\eps}{\varepsilon}
\newcommand{\norm}[1]{\ifthenelse{\equal{#1}{}}{\mbox{$\|\cdot\|$}}{\mbox{$\| #1 \|$}}}
\newcommand{\dual}[2]{\ifthenelse{\equal{#1}{}}{\mbox{$ \left\langle \,\cdot\, , \, \cdot \, \right\rangle $}}{
\mbox{$ \left\langle #1 \,  , \, #2 \right\rangle$}}}

\newcommand{\rg}{\mathrm{rg}}

\newcommand{\CR}{\mathbb{R}}
\newcommand{\CC}{\mathbb{C}}
\newcommand{\CN}{\mathbb{N}}
\newcommand{\one}{\mathbbm{1}}
\newcommand{\weak}{\rightharpoonup}
\newcommand{\minus}{\,\mbox{-}\,}
\renewcommand{\Re}{\mathrm{Re}\,}

\begin{document}
\title{Continuity and Equicontinuity of Transition Semigroups}
\author{Markus Kunze}
\address{Delft Institute of Applied Mathematics, Delft University of Technology, P.O, Box 5031, 2600 GA Delft, The Netherlands}
\email{M.C.Kunze@tudelft.nl}
\thanks{The author was partially supported by the Deutsche Forschungsgesellschaft in the framework of the DFG research training group 1100 at the University of Ulm.}

\begin{abstract}
We study continuity and equicontinuity of semigroups on norming dual pairs
with respect to topologies defined in terms of the duality. In particular, we address the question
whether continuity of a semigroup already implies (local/quasi) equicontinuity. We apply our results
to transition semigroups and show that, under suitable hypothesis on $E$, 
every transition semigroup on $C_b(E)$ which is continuous with
respect to the strict topology $\beta_0$ is automatically quasi-equicontinuous with respect to
that topology. We also give several characterizations of $\beta_0$-continuous semigroups on
$C_b(E)$ and provide a convenient condition for the transition semigroup of a Banach space
valued Markov process to be $\beta_0$-continuous.
\end{abstract}

\maketitle

\section*{Introduction}
An object of central interest in the study of Markov processes is
the transition semigroup of the process. If the Markov process
$(X_t)_{t\geq 0}$ takes values in the measurable space $(E,\Sigma
)$, the state space of the process, then the transition semigroup
$\mathbf{T} = (T(t))_{t\geq 0}$ is a positive contraction semigroup
on the space $B_b(E)$ of all bounded, measurable functions on $E$.
This semigroup contains all information about the transition
probabilities of $X_t$. More precisely, for $t,s \geq 0$ and $A \in
\Sigma$ we have $P (X_{t+s} \in A \, | \, X_s = x ) = (T(t)\one_A)(x)$.\\

Whereas the orbits of the semigroup $\mathbf{T}$ usually bear no
continuity properties, often the restriction of $\mathbf{T}$ to
certain invariant subspaces is continuous in one way or other. The
best known example for this is that of a {\it Feller semigroup}.
Here, $E$ is a locally compact Hausdorff space endowed with the
Borel $\sigma$-algebra.\\
If $C_0(E)$, the space of all continuous
functions on $E$ which vanish at infinity, is invariant under $E$,
then often $\mathbf{T}|_{C_0(E)}$ is strongly continuous. This can
be used to great effect in the study of Markov processes, see
\cite{ito, ek}. If $E$ is not locally compact or if $C_0(E)$ is not
invariant, then one can consider other invariant subspaces. Of
particular interest is the space $C_b(E)$ of all bounded, continuous
functions on $E$. However, even if $C_b(E)$ is invariant under
$\mathbf{T}$, the restriction of $\mathbf{T}$ to $C_b(E)$ is in
general not strongly continuous. But in many cases, see e.g.
\cite{gk01, gvn03, dn96}, the restriction is continuous with
respect to the so-called {\it strict} (or {\it mixed}\/) topology
$\beta_0$, cf. \cite{buck1, buck2, sentilles} and Section 1.3.

For locally compact spaces $E$, {\sc Sentilles} \cite{sentilles2} has studied $\beta_0$-continuous
semigroups on $C_b(E)$ in the framework of equicontinuous semigroups on locally convex spaces \cite{yosida}.
Since $\beta_0$ agrees with the compact-open topology $\tau_{\mathrm{co}}$ on $\norm{}_{\infty}$-bounded
subsets of $C_b(E)$, it is also possible to treat $\beta_0$-continuous semigroups as 
$\tau_{\mathrm{co}}$-continuous semigroups. This point of view was taken by {\sc Cerrai} in \cite{cerrai} and led to the concept of bi-continuous semigroups introduced by 
{\sc K\"uhnemund} in \cite{kuehne}. {\sc Farkas} \cite{farkas09} has used the theory of bi-continuous
semigroups to study transition semigroups on $C_b(E)$, where $E$ is a Polish space, i.e. the 
topology of $E$ is induced by a complete, separable metric. 
It should be noted that transition semigroups on $C_b(E)$ in general do not satisfy the equicontinuity assumption of
\cite{yosida} with respect to $\tau_{\mathrm{co}}$, see \cite[Example 6]{kuehne}. However, in the 
examples in \cite{gk01, gvn03} local equicontinuity with respect to $\beta_0$ holds.\\

In this paper, we study continuity and equicontinuity of semigroups in the framework of semigroups
on norming dual pairs introduced in \cite{k09a}. Thus, in addition to a Banach space $X$, we are given
a closed subspace $Y$ of $X^*$, the norm dual of $X$, which is norming for $X$. We then study semigroups $\mathbf{T}$ on $X$ such that the adjoint semigroup $\mathbf{T}^*$ leaves the space $Y$ invariant. In applications to transition semigroups we will choose $X = C_b(E)$, here the state space $E$ is assumed to be a completely regular
Hausdorff space, and $Y= \mathcal{M}_0(E)$, the space of all bounded Radon measures on $E$. In this context the assumption that $\mathbf{T}^*\mathcal{M}_0(E) \subset \mathcal{M}_0(E)$ is quite natural and has a stochastic interpretation. Namely, if $\mathbf{T}$ is the transition semigroup of the Markov process $(X_t)_{t\geq 0}$ and we put $\mathbf{T}' = \mathbf{T}^*|_{\mathcal{M}_0(E)}$, then $\mathbf{T}'$ gives the distribution of the random elements $X_t$, i.e. if $X_s\sim \mu \in \mathcal{M}_0(E)$ then $X_{t+s}\sim T(t)'\mu$.

In Sections 2 and 3 we will work on general norming dual pairs and study continuity and equicontinuity
with respect to general locally convex topologies defined in terms of the duality, see Section 1.1.
This generality allows us to consider continuity with respect to various topologies. In particular,
if we choose $\tau = \norm{}$, then we obtain strongly continuous semigroups. On the norming dual pair
$(C_b(E),\mathcal{M}_0(E))$ not only the strict topology $\beta_0$ but also the weak topology $\sigma (C_b(E),\mathcal{M}_0(E))$ is of interest. This topology is connected with the concept of
{\it bounded and pointwise convergence}, see \cite[Section 3.4]{ek}. {\sc Priola} \cite{pri} has used
this continuity concept to study transition semigroups.

If we additionally impose certain equicontinuity assumptions, then it not surprising that we can prove a generation theorem for such semigroups. The more interesting question is whether equicontinuity assumptions are restrictive or whether, at least for certain
topologies, these assumptions are satisfied automatically. We will address this question in Section 3
and give some abstract examples where this is the case. In Section 4 we apply our results to 
transition semigroups. We will prove that if $E$ is a Polish space, then every $\beta_0$-continuous semigroup on $(C_b(E),\mathcal{M}_0(E))$ is locally $\beta_0$-equicontinuous. 
A variant of this result has been
obtained independently by {\sc Farkas} \cite{farkas09}. However, we will also prove that this 
result remains valid for {\it positive} semigroups, whenever $E$ is a
so-called T-space (see the definition in Section 4). In the main result of Section 4,
Theorem \ref{t.cbmucont}, we give various equivalent conditions for
a semigroup on $(C_b(E), \mathcal{M}_0(E))$ to be $\beta_0$-continuous. In the concluding
Section 5 we discuss several examples and give a convenient condition for the transition 
semigroup of a Banach space valued Markov process to be $\beta_0$-continuous.

\section{Preliminaries and Notation}

\subsection{Dual pairs} Throughout this paper we will be working on
dual pairs and use locally convex topologies defined in terms of the
duality. We briefly recall some results from the theory and fix some
notation. Our main reference are Chapters 20 and 21 of
\cite{koethe}. A {\it dual pair} is a triple $(X, Y, \dual{}{})$
where $X$ and $Y$ are vector spaces over the same field $\mathbb{K}
= \CR$ or $\CC$ and $\dual{}{}$ is a bilinear form from $X\times Y$
to $\mathbb{K}$ which separates points, i.e. $\dual{x}{y} = 0$ for
all $x \in X$ implies $y = 0$ and $\dual{x}{y} = 0$ for all $y \in
Y$ implies $x = 0$. We may define locally convex topologies on
$X$ as follows. If $M \subset Y$ is bounded, i.e. $\sup_{y\in
M}|\dual{x}{y}| < \infty$ for all $x\in X$, then $p_M(x) :=
\sup_{y\in M} |\dual{x}{y}|$ defines a seminorm on $X$. If
$\mathfrak{M}$ is a collection of bounded subsets of $Y$, then the
collection of seminorms $(p_M)_{M\in\mathfrak{M}}$ defines a locally
convex topology on $X$ if and only if for every $x \in X$ there
exists some $M \in \mathfrak{M}$ such that $p_M(x) \neq 0$ (we say that 
$\mathfrak{M}$ is {\it separating}). If $\mathfrak{M}$ is a separating
collection of bounded subsets of $Y$, then $\tau_{\mathfrak{M}}$ denotes the locally convex
topology induced by the seminorms $(p_M)_{M\in\mathfrak{M}}$.\\

A locally convex topology $\tau$ on $X$ is called {\it consistent}
if $(X,\tau )' =Y$, i.e. every $\tau$-continuous linear functional
$\varphi$ on $X$ is of the form $\varphi (x) = \dual{x}{y}$ for some
$y \in Y$. By the Mackey-Arens theorem, \cite[21.4 (2)]{koethe},
every consistent topology is of the form $\tau_{\mathfrak{M}}$ for a
suitable collection $\mathfrak{M}$. Furthermore, there exists a
coarsest consistent topology, namely the {\it weak topology} $\sigma
(X,Y) = \tau_{\mathfrak{F}}$, where $\mathfrak{F}$ denotes the
collection of all finite subsets of $Y$, and a finest consistent
topology, namely the {\it Mackey topology} $\mu (X,Y) =
\tau_{\mathfrak{K}}$, where $\mathfrak{K}$ denotes the collection of
all absolutely convex, $\sigma (Y,X)$-compact subsets of $Y$. We note
that every topology $\tau_{\mathfrak{M}}$ is finer than the weak
topology $\sigma (X,Y)$. To simplify notation, we will write
$\sigma$ (resp. $\mu$) for $\sigma (X,Y)$ (resp. $\mu (X,Y)$) and denote $\sigma$-convergence on 
$X$ by $\weak$. We will write $\sigma'$ (resp. $\mu'$) for $\sigma
(Y,X)$ (resp. $\mu(Y,X)$) and denote $\sigma'$-convergence on $Y$ by $\weak'$.

\subsection{Norming dual pairs}
\begin{defn}
A {\it norming dual pair} is a dual pair $(X,Y, \dual{}{})$ where
$X$ and $Y$ are Banach spaces and we have $\norm{x} = \sup\{\,
|\dual{x}{y}|\, : \, y \in Y\, , \, \norm{y}\leq 1\,\}$ and
$\norm{y} = \sup\{\, |\dual{x}{y}|\, : \, x \in X\, , \,
\norm{x}\leq 1\,\}$.
\end{defn}
In what follows, we will often write $(X,Y)$ instead of
$(X,Y,\dual{}{})$ if the duality pairing is understood. It is easy
to see that if $(X,Y)$ is a norming dual pair, then $Y$ is
isometrically isomorphic to a closed subspace of $X^*$, the norm dual of $X$.
We will often identify $Y$ with this closed subspace of $X^*$.

It is an easy but crucial consequence of the definition that if $(X,Y)$ is a
norming dual pair, then on $X$ and $Y$ the notions of weak (i.e. $\sigma$- or
$\sigma'$-) boundedness and of norm boundedness coincide, cf.
\cite{k09a}. It follows that the norm topology on $X$ is equal to $\tau_{\mathfrak{B}}$, 
where $\mathfrak{B}$ denotes the collection of all bounded subsets of $Y$. 
The norm topology is in general not consistent, 
but it is finer than any topology $\tau_{\mathfrak{M}}$.
In particular, if $T$ is a $\tau_{\mathfrak{M}}$-continuous linear
operator, then it is $\norm{}$-continuous. It is proved in
\cite{k09a} that a $\norm{}$-continuous linear operator $T$ is
$\sigma$-continuous if and only if its norm-adjoint $T^*$ leaves the
space $Y$ invariant. By \cite[21.4 (6)]{koethe}, a linear operator
is $\sigma$-continuous if and only if it is $\mu$-continuous. If
$\tau$ is a consistent locally convex topology, then every
$\tau$-continuous linear operator is $\sigma$-continuous. The
converse is not true in general. If $\tau$ is any (not necessarily
consistent) locally convex topology on $X$, we write $L(X,\tau)$ for
the algebra of $\tau$-continuous linear operators on $X$. For $\tau
=\norm{}$ we merely write $L(X)$ instead of $L(X,\norm{})$. If $T
\in L(X,\sigma )$, we write $T^*$ for its norm adjoint and $T'$ for
its $\sigma$-adjoint. Note that $T'=T^*|_Y$.

\subsection{The dual pair $\mathbf{(C_b(E),\mathcal{M}_0(E))}$}
Our main example for applications is the norming dual pair $(C_b(E),
\mathcal{M}_0(E))$. Here, $E$ is a completely regular Hausdorff
space and $C_b(E)$ denotes the Banach space of all bounded,
continuous functions form $E$ to $\CC$ endowed with the supremum
norm. A positive measure $\mu$, defined on the Borel
$\sigma$-algebra $\mathcal{B}(E)$, is called a {\it Radon measure}
if for all $A \in \mathcal{B}(E)$, we have $\mu (A) = \sup\{ \mu (K)
\, : \, K \subset A\, , K \,\, \mbox{compact} \, \}$. If $\mu$ is a
complex measure on $\mathcal{B}(E)$, then its total variation $|\mu
|$ is defined by $|\mu |(A) =
\sup_{\mathcal{Z}}\sum_{B\in\mathcal{Z}}|\mu (B)|$, where the
supremum is taken over all finite partitions $\mathcal{Z}$ of $A$
into pairwise disjoint measureable sets. A complex measure $\mu$ is
called a Radon measure, if $|\mu|$ is a Radon measure. Note that if
$E$ is a Polish space, then every measure on $\mathcal{B}(E)$ is a Radon measure.
$\mathcal{M}_0(E)$ denotes the Banach space of all bounded Radon measures
on $E$, endowed with the total variation norm $\norm{\mu} := |\mu
|(E)$. It is proved in \cite{k09a} that $(C_b(E),
\mathcal{M}_0(E))$ is a norming dual pair with respect to the
duality $\dual{f}{\mu} = \int_E f \, d\mu $. If $T \in L(C_b(E),
\sigma )$, then $T$ has the representation $Tf(x) = \int_E f(y)\,
k(x, dy)$. Here, $k(x, \cdot ) = T'\delta_x$, where $\delta_x$
denotes the Dirac measure in $x$. We will call $k$ the {\it kernel
associated with} $T$. The question whether $k(\cdot , A)$ is
measurable for all $A \in \mathcal{B}(E)$ is
discussed in \cite{k09a}.\\

The strict topology $\beta_0$ on $C_b(E)$ is defined as follows:\\
Denote by $\mathcal{F}_0(E)$ the space of all bounded functions on
$E$ which vanish at infinity, i.e. given $\eps > 0$, we find a compact set $K
\subset E$ such that $|f(x)|\leq \eps$ for all $x\not\in K$. The
{\it strict topology} $\beta_0$ on $C_b(E)$ is the locally convex
topology generated by the set of seminorms $(p_{\varphi})_{\varphi
\in \mathcal{F}_0(E)}$, where $p_{\varphi}(f) := \norm{\varphi f}_{\infty}$.

This definition is taken from \cite{jarchow}. It generalizes the
definition given by {\sc Buck} \cite{buck1, buck2} for locally compact spaces $E$. 
By \cite[Theorem 7.6.3]{jarchow}, $(C_b(E), \beta_0 )' = \mathcal{M}_0(E)$, i.e.
$\beta_0$ is a consistent topology. Furthermore, $(C_b(E),\beta_0 )$
is complete if and only if $C(E)$, the space of all continuous
functions on $E$, is complete with respect to $\tau_{\mathrm{co}}$,
see Theorems 4 and 9 in Section 3.6 of \cite{jarchow}. In particular, if
$E$ is a metric space or a locally compact space, then $(C_0(E), \beta_0)$
is complete. {\sc Sentilles} \cite{sentilles} has considered several strict topologies
yielding different spaces of measures as dual spaces. We will recall
some results from \cite{sentilles} in Section 4.

\section{Semigroups and their Generators}

We now study semigroups on norming dual pairs. As a matter of fact, several interesting
properties of such semigroups can be proved without continuity assumptions, merely imposing
integrability assumptions. This leads to the concept of integrable semigroups on norming dual
pairs. Such semigroups are studied in \cite{k09a} and we content ourselves with recalling the
definition and collecting some of the results from \cite{k09a} in Propositions \ref{p.pseudoresolvent}
and \ref{t.rangecharacterisation} below.

\begin{defn}\label{d.semigroup}
Let $(X,Y)$ be a norming dual pair. A {\it semigroup} on $(X,Y)$ is a
family $\mathbf{T} = (T(t))_{t\geq 0} \subset L(X, \sigma )$ such that
\begin{enumerate}
 \item $\mathbf{T}$ is a {\it semigroup}, i.e. $T(0)=0$ and $T(t+s) = T(t)T(s)$ for all $t,s \geq 0$.
 \item $\mathbf{T}$ is {\it exponentially bounded}, i.e. there exist $M\geq 1$ and $\omega \in \CR$ such that
 $\norm{T(t)}\leq Me^{\omega t}$ for all $t\geq 0$. We then say that $\mathbf{T}$ is of
{\it type} $(M,\omega )$.
\end{enumerate}
A semigroup $\mathbf{T}$ of type $(M,\omega )$ is called {\it integrable} if
\begin{itemize}
 \item[(3)] for all $\lambda$ with $\Re\lambda > \omega$, there exists an operator $R(\lambda ) \in L(X,\sigma )$
such that
\begin{equation}\label{eq.resolvent}
 \dual{R(\lambda )x}{y} = \int_0^{\infty} e^{-\lambda t}\dual{T(t)x}{y} \,\, .
\end{equation}
In particular, we assume that all the integrals on the right hand side exist. $\mathbf{R}
= (R(\lambda ))_{\Re\lambda > \omega }$ is called the {\it Laplace transform of} $\mathbf{T}$.
\end{itemize}
\end{defn}

\begin{prop}\label{p.pseudoresolvent}
Let $\mathbf{T}$ be an integrable semigroup of type $(M,\omega )$ with Laplace transform $\mathbf{R}$.
\begin{enumerate}
 \item $\mathbf{R}$ is a pseudoresolvent and every $R(\lambda )$ commutes with every $T(t)$.
 \item We have $\norm{(\Re \lambda - \omega )^kR(\lambda )^k} \leq M$ for all $\Re\lambda > \omega$
and $k \in \CN$.
 \item If $\rg \, \mathbf{R}$ is $\sigma$-dense in $X$, then $\mathbf{R}$ determines 
$\mathbf{T}$ uniquely.
More precisely, if $\tilde{\mathbf{T}}$ is an integrable semigroup on $(X,Y)$ having the same Laplace
transform $\mathbf{R}$, then $\mathbf{T} =\tilde{\mathbf{T}}$.
\end{enumerate}
\end{prop}

Recall that a pseudoresolvent is a map $R$ from some nonempty set
$\Omega \subset \CC$ to $L(X,\norm{})$, such that $R(\lambda ) -
R(\mu ) = (\mu -\lambda )R(\lambda )R(\mu )$ for all $\lambda , \mu
\in \Omega$. It is well known that for a given pseudoresolvent
$(R(\lambda ))_{\lambda \in \Omega}$ there exists a unique
multivalued operator $\mathcal{A}$ such that $R(\lambda ) = (\lambda
- \mathcal{A})^{-1}$ for all $\lambda \in \Omega$. In particular,
the range $\rg R(\lambda )$ and
the kernel $\ker R(\lambda )$ are independent of $\lambda \in \Omega$.

The following proposition gives a characterization of the operator $\mathcal{A}$. The integrals
appearing are to be understood as $Y$-integrals, see \cite{k09a}. More precisely, if
$f: I \to X$ is a function defined on some interval $I$ such that $\dual{f(\cdot )}{y}$ is integrable
for every $y \in Y$, then the $Y$-integral $\int_I f(t) dt$ denotes the unique element $\varphi \in Y^*$ such that $\varphi (y) = \int_I \dual{f(t)}{y} \, dt$ for all $y \in Y$.
In the proposition below, $\int_I f(t) \, dt$ will actually be an
element of $X$ which is considered as a closed subspace of $Y^*$.
Hence,  $\int_0^tf(t)\, dt =x \in X$  if and only if $\dual{x}{y} =
\int_I\dual{f(t)}{y}\, dt$ for all $y \in Y$. However, even if
$\int_If(t)\, dt \in X$, the integral in general does not exist as
a Bochner or as a Pettis integral.

\begin{prop}\label{t.rangecharacterisation}
Let $\mathbf{T}$ be an integrable semigroup on the norming dual pair
$(X,Y)$ with Laplace transform $R(\lambda ) = (\lambda -
\mathcal{A})^{-1}$.
\begin{enumerate}
\item The following are equivalent.
\begin{enumerate}
\item $x \in D(\mathcal{A})$ and $z \in \mathcal{A}x$;
\item for every $t> 0$ and $y \in Y$ we have
\begin{equation}\label{eq.awf}
 \int_0^tT(s)z\, ds = T(t)x - x \,\, .
\end{equation}
\end{enumerate}
\item For $x \in X$ and $t>0$ we have $\int_0^t T(s)x \, ds \in D(\mathcal{A})$ and
\[ T(t)x -x \in \mathcal{A}\int_0^t T(s)x \, ds \,\, . \]
\end{enumerate}
\end{prop}

\begin{rem}
 It follows from \eqref{eq.awf} that $t \mapsto T(t)x$ is $\norm{}$-continuous for
every $x \in D(\mathcal{A})$.  Indeed, if $x\in D(\mathcal{A})$ and $z \in \mathcal{A}x$,
then for every $t_0>0$ we have $C:=\sup_{t\leq t_0}\norm{T(t)} < \infty$.
Hence, \eqref{eq.awf} implies
that
\[
 |\dual{T(t)x -T(s)x}{y}| \leq \int_s^t |\dual{T(r)z}{y}|\, dr \leq  |t-s| 2C \norm{z}\cdot
  \norm{y} \,\, .
 \]
for $t,s \leq t_0$ and $y\in Y$. Taking the supremum over $y\in Y$ with $\norm{y}\leq 1$,
$\norm{}$-continuity
of $t\mapsto T(t)x$ follows.
\end{rem}

\begin{defn}\label{d.generator}
Let $\mathbf{T}$ be an integrable semigroup on the norming
dual pair $(X,Y)$ such that the Laplace transform $\mathbf{R}$ of $\mathbf{T}$ is injective.
Then the unique (single valued) operator $A$ such that $R(\lambda ) = (\lambda - A)^{-1}$
is called the {\it generator} of $\mathbf{T}$. In this case we say that $\mathbf{T}$ {\it has
a generator} or that $\mathbf{T}$ is a {\it semigroup with generator} ($A$).
\end{defn}

If $\tau$ is a locally convex topology on $X$, then, as usual, 
an operator $A$ on $X$ is called {\it $\tau$-closed} if the graph of $A$ is closed in 
$X\times X$ with respect to $\tau\times \tau$. 
If $\tau$ is a consistent topology, then, by the Hahn-Banach theorem, an operator
$A$ is $\tau$-closed if and only if it is $\sigma$-closed.
Furthermore, a $\sigma$-closed operator is automatically norm closed.  
For an operator $A$, we denote its resolvent set by $\rho (A)$ and for $\lambda \in \rho (A)$
we write $R(\lambda , A)$ for the resolvent of $A$ in $\lambda$. We define
\[ \rho_{\sigma }(A) := \{ \lambda \in \rho (A) \, : \,
R(\lambda , A)\in L(X,\sigma ) \, \}\,\, . \] 
It is an open question whether $\rho_{\sigma}(A) = \rho (A)$ for a $\sigma$-closed
operator $A$. For a $\sigma$-densely defined, $\sigma$-closed
operator, the $\sigma$-adjoint of $A$ is denoted by $A'$.

\begin{prop}\label{p.generator}
Let $\mathbf{T}$ be an integrable semigroup of type $(M,\omega )$ with
generator $A$. Then $A$ is a $\sigma$-closed operator with $\{ \Re\lambda > \omega \}
\subset \rho_{\sigma}(A)$. Furthermore, for $\Re\lambda > \omega$ and $k\in \CN_0$
we have
\begin{equation}\label{eq.hilleyosida}
\norm{R(\lambda , A)^k} \leq\frac{M}{(\Re\lambda - \omega )^k} \,\, .
\end{equation}

The operator $A$ is $\sigma$-densely defined if and only if $\mathbf{T}'$ is a semigroup
with generator.
\end{prop}

\begin{proof}
Since the resolvent of $A$ is the Laplace transform of $\mathbf{T}$
and since the Laplace transform consists of $\sigma$-continuous
operators, $\{ \Re\lambda \, : \, \lambda > \omega \} \subset
\rho_{\sigma}(A)$. In particular, $A$ is $\sigma$-closed. Estimate
\eqref{eq.hilleyosida} follows from Proposition
\ref{p.pseudoresolvent}. Now assume that $A$ is $\sigma$-densely
defined. In this case, the $\sigma$-adjoint $A'$ of $A$ is well
defined and $R(\lambda , A)' = R(\lambda , A')$, as is easy to see.
Since clearly $\dual{x}{R(\lambda , A)'y} = \int_0^{\infty}e^{-\lambda t}\dual{x}{T(t)'y}\, dt$
for all $x \in X$ and $y \in Y$, it follows that
$A'$ is the generator of $\mathbf{T}'$. $\quad$ Conversely assume
that $\mathbf{T}'$ has a generator $B$. As the Laplace transform of
$\mathbf{T}'$ is $R(\lambda , A)'$, we find $R(\lambda , A)'
= R(\lambda , B)$. If $y\in Y$ vanishes on $D(A)$, then $0 =
\dual{R(\lambda , A)x}{y} = \dual{x}{R(\lambda , B)y}$ for all $x
\in X$. It follows that $R(\lambda , B)y = 0$. As $R(\lambda , B)$
is injective by hypothesis, $y = 0$ follows. By the Hahn-Banach
theorem, $D(A)$ is $\sigma$-dense in $X$.
\end{proof}

We now turn to continuous semigroups.

\begin{defn}
Let $\mathbf{T}$ be a semigroup on $(X,Y)$ and $\tau$ be a locally convex
topology on $X$. Then $\mathbf{T}$ is called {\it $\tau$-continuous (at 0)} if
for all $x \in X$ the map $t \mapsto T(t)x$ is $\tau$-continuous (at 0).
\end{defn}

\begin{rem}\label{rem1}
\begin{enumerate}
\item Using the uniform boundedness principle, it can be shown that if
$\mathbf{T}$ is a semigroup on $(X,Y)$ which is $\sigma$-continuous
at $0$, then $\mathbf{T}$ is automatically exponentially bounded, i.e. condition
(2) in Definition \ref{d.semigroup} is automatically satisfied, see \cite{k09a}.
\item In the remainder of this section, we will assume integrability of semigroups, 
i.e. condition (3) in Definition \ref{d.semigroup}, also under continuity assumptions.
This is due to the fact that $\sigma$-continuity at $0$ in general does not imply
integrability of a semigroup, see the example in Section 5.2. However, it is proved
in \cite{k09a} that if $E$ is a complete metric space, then every semigroup on
$(C_b(E),\mathcal{M}_0(E))$ which is $\sigma$-continuous at $0$ is integrable.
\end{enumerate}
\end{rem}

Our definition of the generator via the Laplace transform is in the spirit of \cite{abhn}.
The following Theorem shows that, under continuity assumptions, this ``integral definition'' 
coincides with the ``differential definition'' of the generator, see e.g. \cite{en}. 

\begin{thm}\label{t.resolventwc}
Let $\mathbf{T}$ be an integrable semigroup on $(X,Y)$ of type
$(M,\omega )$ and $\mathfrak{M}$ be a separating collection of bounded subsets of $Y$.
If $\mathbf{T}$ is $\tau_{\mathfrak{M}}$-continuous at 0, then
$\tau_{\mathfrak{M}}\minus\lim_{\lambda \to \infty}\lambda R(\lambda
)x = x$. In particular, the Laplace transform of $\mathbf{T}$ is
injective and $\mathbf{T}$ has a generator $A$ such that $D(A)$ is
sequentially $\tau_{\mathfrak{M}}$-dense in $X$. Furthermore, the
generator is given by
\[
 D(A) = \left\{ x \in X \, : \tau_{\mathfrak{M}}\minus\lim_{h\downarrow 0}\Delta_hx\,\,\,
\mbox{exists}\,\right\}\quad , \quad Ax = \tau_{\mathfrak{M}}\minus\lim_{h\downarrow 0}\Delta_hx \,\,,
\]
where $\Delta_hx := h^{-1}(T(h)x - x)$.
\end{thm}

\begin{proof}
Let $x \in X$, $S \in \mathfrak{M}$ and $\eps > 0$ be given.
Since $S$ is bounded, there exists $C > 0$ such that $\norm{y} \leq
C$ for all $y \in S$.
By $\tau_{\mathfrak{M}}$-continuity at 0, there
exists $t_0 > 0$ such that $|\dual{T(t)x -x}{y}| \leq \eps$ for all $t
\leq t_0$ and  $y \in S$.
Now for $\lambda > \max\{\omega , 0\}$ and $y \in S$ we have
\begin{eqnarray*}
\sup_{y\in S} |\dual{\lambda R(\lambda )x - x}{y}| & = &
\sup_{y\in S}\left|\int_0^{\infty}\dual{\lambda e^{-\lambda t}T(t)x - \lambda
e^{-\lambda t}x}{y}\, dy\right|\\
& \leq & \sup_{y\in S}\int_0^{t_0}\lambda e^{-\lambda t}|\dual{T(t)x-x}{y}|\, dt\\
&& + \int_{t_0}^{\infty} \lambda e^{-\lambda t}(1+ M e^{\omega t})
C\cdot \norm{x}\, dt\\
& \leq & \eps\left( 1-e^{-\lambda t_0}\right) +
C\cdot\norm{x}\left(e^{-\lambda t_0}
+ \frac{\lambda\cdot M}{\lambda -\omega}e^{(\omega -\lambda )t_0}\right)\\
& \to &\eps
\end{eqnarray*}
as $\lambda \to \infty$. Since $S \in \mathfrak{M}$ was arbitrary, the first part is proved.
Now denote the generator of $\mathbf{T}$ (in the sense of Definition \ref{d.generator}) by
$B$ and let $A$ be the operator in the statement. If $x \in D(B)$, then, by Proposition \ref{t.rangecharacterisation}, we have
\[ |\dual{\Delta_hx -Bx}{y} | \leq \frac{1}{h} \int_0^h |\dual{T(s)Bx -Bx}{y}| \, ds \,\, ,\]
for every $y \in Y$. Now let $S \in \mathfrak{M}$ and $\eps > 0$ be given. 
Choose $t_0>0$ such that $p_S(T(s)Bx -Bx) \leq \eps$, for all $0\leq s \leq t_0$. 
Then,
\[ p_S\left(\Delta_hx - Bx \right) \leq  \frac{1}{h}\int_0^h \eps \, ds = \eps \,\, , \]
for all $0 \leq s \leq t_0$.
This proves that $x \in D(A)$ and that $Ax = Bx$. Conversely
suppose that $x \in D(A)$. Since $\tau_{\mathfrak{M}}$ is finer than
$\sigma$ it follows that $\Delta_hx \weak Ax$ as $h \downarrow 0$.
Since every operator $T(s)$ is $\sigma$-continuous, $T(s)\Delta_hx
\weak T(s)Ax$ for every $s \geq 0$. Furthermore, $(\Delta_hx)_{h
\leq 1}$ is norm bounded. Indeed, for every $y \in Y$, the set
$\{\dual{\Delta_hx}{y}\, : \, h\leq 1 \}$ is bounded. Hence, by the
uniform boundedness principle, $\sup_{h\leq 1}\norm{\Delta_hx}_{Y^*}
< \infty$. However, since $X$ embeds isometrically into
$Y^*$, we have $\norm{\Delta_h x}_{Y^*} = \norm{\Delta_h x}_X$ for
every $h > 0$.\\
Now fix $t>0$ and $y \in Y$. Put $I_t := \int_0^t T(s)x \, ds$. Then
\begin{eqnarray*}
 \int_0^t\dual{T(s)Ax}{y}\, ds & = & \lim_{h\downarrow 0}\int_0^t \dual{T(s)\Delta_hx}{y}\, ds
 =  \lim_{h\downarrow 0} \dual{\Delta_h I_t}{y} \,\, ,
\end{eqnarray*}
by the dominated convergence theorem. Note that $I_t \in D(B)$ and
$\dual{BI_t}{y} = \dual{T(t)x -x}{y}$ by Proposition \ref{t.rangecharacterisation}. 
Since $B\subset A$, it follows
that
\[ \int_0^t\dual{T(s)Ax}{y} = \lim_{h\downarrow 0 } \dual{\Delta I_t}{y} 
 = \dual{BI_t}{y} = \dual{T(t)x-x}{y} \,\, .
\]
Thus Proposition \ref{t.rangecharacterisation} implies that $x \in D(B)$ and $Bx = Ax$.
\end{proof}

\begin{rem}
Assume in addition to the hypothesis of Theorem \ref{t.resolventwc} that the semigroup $\mathbf{T}$
is $\tau_{\mathfrak{M}}$-continuous. Then, arguing similar as in the proof of Theorem \ref{t.resolventwc}, it is easy to see that $x \in D(A)$ if and only if $t \mapsto T(t)x$ is
$\tau_{\mathfrak{M}}$-differentiable. In this case we have $\frac{d}{dt}T(t)x = T(t)Ax$. 
Note however that $\tau_{\mathfrak{M}}$-continuity at 0 does not imply 
$\tau_{\mathfrak{M}}$-continuity.
\end{rem}

We now give a characterization of continuous semigroups.

\begin{prop}\label{p.weaklycontinuous}
Let $(X,Y)$ be a norming dual pair and $\mathfrak{M}$ be a separating collection
of bounded subsets of $Y$. Furthermore, let $\mathbf{T}$ be an integrable semigroup on $(X,Y)$ with
generator. Then the following are equivalent:
\begin{enumerate}
\item $\mathbf{T}$ is $\tau_{\mathfrak{M}}$-continuous;
\item For all $t_0 > 0$ and every $x \in X$ the set
$\{ T(t)x \, : \, t \in [0,t_0] \}$
is $\tau_{\mathfrak{M}}$-compact.
\item For some $t_0>0$ and every $x \in X$ the set
$\{ T(t)x\, : \, t \in [0,t_0] \}$ is relatively countably $\tau_{\mathfrak{M}}$-compact.
\end{enumerate}
\end{prop}

\begin{proof}
(1) $\Rightarrow$ (2) and (2) $\Rightarrow$ (3) are trivial. For (3) $\Rightarrow$ (1) 
suppose that $t \mapsto T(t)x$ is not $\tau_{\mathfrak{M}}$-continuous at 
$ t \in [0,t_0]$. Then there exists a $\tau_{\mathfrak{M}}$-continuous seminorm $p$, an
$\eps >0$ and a sequence $(t_n) \subset [0,t_0]$ converging to $t$ such that $p(T(t_n)x-T(t)x)\geq \eps$
for all $n \in \CN$. By hypothesis, the sequence $T(t_n)x$ has an accumulation point $z \in X$. Thus there exists a subnet $t_\alpha$ of $t_n$ such that $T(t_\alpha)x \stackrel{\tau_{\mathfrak{M}}}{\to} z$. Since $\tau_{\mathfrak{M}}$ is finer than $\sigma$ we have $T(t_{\alpha})x \weak z$
As $R(\lambda )$ is $\sigma$-continuous and commutes with the semigroup, we have
\[ R(\lambda )z = \sigma\minus\lim R(\lambda )T(t_\alpha)x =
\sigma\minus\lim T(t_\alpha)R(\lambda )x = T(t)R(\lambda )x = R(\lambda )T(t)x \,\, ,\]
as $s\mapsto T(s)R(\lambda )x$ is $\norm{}$-continuous. Since $R(\lambda )$ is injective, it follows
that $z = T(t)x$. But then $p(T(t_\alpha )x -T(t)x) \to 0$, a contradiction.
This proves that $t \mapsto T(t)x$ is $\tau_{\mathfrak{M}}$-continuous on $[0,t_0]$. 
Using the semigroup law and that $\{T(t)T(t_0)x \, : \, t \in [0,t_0] \}$ is relatively countably
compact, it follows that $t\mapsto T(t)x$ is $\tau_{\mathfrak{M}}$-continuous on $[0, 2t_0]$.
Inductively we obtain continuity for all times.
\end{proof}

\section{Equicontinuity}

In the context of semigroups, several types of equicontinuity assumptions have been discussed
in the literature. We briefly recall the definition.

\begin{defn}
Let $(X, \tau )$ be a locally convex space. A set
$\mathcal{S}\subset L(X, \tau )$ is called {\it equicontinuous}, if
for every $\tau$-continuous seminorm $p$, there exists a
$\tau$-continuous seminorm $q$ such that $p(Tx) \leq q(x)$ for all
$x \in X$ and $T \in \mathcal{S}$. 
A semigroup $\mathbf{T}$ of $\tau$-continuous operators
is called {\it locally $\tau$-equicontinuous}, if $\{ T(t) \, : \, t
\in [0,t_0] \}$ is $\tau$-equicontinuous for all $t_0 > 0$. It is
called {\it (globally) $\tau$-equicontinuous}, if $\{ T(t) \, : \, t
\geq 0 \}$ is $\tau$-equicontinuous. If for some $\alpha \in \CR$
the rescaled semigroup $\mathbf{T}_\alpha := (e^{-\alpha t}T(t))_{t\geq 0}$
is $\tau$-equicontinuous, then $\mathbf{T}$ is called {\it ${}_{\alpha}$quasi-$\tau$-equicontinuous}.
We will say that $\mathbf{T}$ is {\it quasi-$\tau$-equicontinuous}, if it is
${}_\alpha$quasi-$\tau$-equicontinuous for some $\alpha \in \CR$.
\end{defn}

Obviously, every quasi-$\tau$-equicontinuous semigroup (in particular, every $\tau$-equicontinuous
semigroup) is locally equicontinuous. The converse is not true in general. 

\begin{exmp}
Consider the norming dual pair $(C_b(\CR), \mathcal{M}_0(\CR ))$. The compact open
topology $\tau_{\mathrm{co}}$ is of the form $\tau_{\mathfrak{M}}$. More precisely,
$\mathfrak{M}$ is the separating collection of sets of the form $\{\delta_x \,:\, x \in K\}$,
where $\delta_x$ denotes the Dirac measure in $x$ and $K$ is a compact subset of $\CR$.
The shift semigroup $\mathbf{T}$, defined by $T(t)f(x) = f(x+t )$ is locally $\tau_{\mathrm{co}}$-equicontinuous but not quasi-$\tau_{\mathrm{co}}$-equicontinuous.
\end{exmp}

\begin{prop}\label{p.taucont}
Let $\mathbf{T}$ be an integrable semigroup on $(X,Y)$ with
generator $A$ and $\mathfrak{M}$ be a separating collection of bounded subsets of $Y$.
If $\mathbf{T}$ is locally $\tau_{\mathfrak{M}}$-equicontinuous and $D(A)$ is
$\tau_{\mathfrak{M}}$-dense in $X$, then $\mathbf{T}$ is
$\tau_{\mathfrak{M}}$-continuous. 
\end{prop}

\begin{proof}
We first prove that $X_0 := \{x \in X \, : \, t\mapsto T(t)x
\,\,\,\mbox{is}\,\, \tau\minus\mbox{continuous} \, \}$ is
$\tau_{\mathfrak{M}}$-closed in $X$. Let $x$ be an accumulation
point of $X_0$, $t_0 \geq 0$ and $p$ be a
$\tau_{\mathfrak{M}}$-continuous seminorm. Pick a
$\tau_{\mathfrak{M}}$-continuous seminorm $q$ such that $p(T(t)z)
\leq q(z)$ for all $t \in [0, t_0+1]$ and $z \in X$. Given $\eps
>0$, we find $x_0 \in X_0$ such that $q(x-x_0) \leq \eps$.
Since $x_0 \in X_0$, there exists $0<\delta < 1$ such that
$p(T(t_0)x_0 - T(t)x_0) \leq \eps$ for all $|t-t_0| \leq \delta$.
Now
\begin{eqnarray*}
p(T(t_0)x - T(t)x) & \leq & p(T(t_0)x - T(t_0)x_0) + p(T(t_0)x_0 -T(t)x_0)\\
&& + p(T(t)x_0 - T(t)x )\\
& \leq & 2q(x-x_0) + p(T(t_0)x_0 - T(t)x_0) \leq 3\eps \,\, ,
\end{eqnarray*}
for all $|t-t_0| \leq \delta$. This proves that $x\in X_0$, whence $X_0$ is $\tau_{\mathfrak{M}}$-closed.
For $x \in D(A)$, $t\mapsto T(t)x$ is $\norm{}$-continuous and hence
$\tau_{\mathfrak{M}}$-continuous. Thus $D(A) \subset X_0$. As $D(A)$ is
$\tau_{\mathfrak{M}}$-dense, $X_0 = X$ follows.
\end{proof}

For $\tau_{\mathfrak{M}}= \norm{}$, we note that local norm-equicontinuity of 
a semigroup
is equivalent with exponential boundedness. Hence from Theorem
\ref{t.resolventwc} and Proposition \ref{p.taucont} we immediately obtain the
following characterization.

\begin{cor}\label{c.weakequalstrong}
Let $\mathbf{T}$ be an integrable semigroup on $(X,Y)$. The
following are equivalent:
\begin{enumerate}
\item $\mathbf{T}$ is strongly continuous;
\item $\mathbf{T}$ has a $\norm{}$-densely defined generator.
\end{enumerate}
\end{cor}

For quasi-equicontinuous semigroups, we obtain the following generation result.

\begin{thm}\label{t.hilleyosida}
Let $(X,Y)$ be a norming dual pair, $\tau$ be a
consistent topology on $X$ which is sequentially complete, and $A$
be a $\sigma$-closed operator on $X$. The following are equivalent.
\begin{enumerate}
\item $A$ is the generator of a $\tau$-continuous, ${}_\alpha$quasi-$\tau$-equicontinuous, integrable 
semigroup $\mathbf{T}$ on $(X,Y)$ of type $(M,\omega )$; 
\item $A$ is a sequentially $\tau$-densely defined operator
such that
\begin{enumerate}
\item $\{ \lambda \in \CR \, : \, \lambda > \omega \} \subset \rho_{\sigma}(A)$ and
\[ \norm{(\lambda - \omega )^kR(\lambda, A)^k} \leq M \quad
\forall \, \lambda > \omega \,\, , k \in \CN \]
\item The set
\[ \left\{ (\lambda - \alpha )^k R(\lambda , A)^k \, : \,
\lambda > \alpha \, , \, k\in \CN \, \right\} \] is
$\tau$-equicontinuous.
\end{enumerate}
\end{enumerate}
\end{thm}

\begin{proof}
(1) $\Rightarrow$ (2): $A$ is sequentially
$\tau$-densely defined by Theorem
\ref{t.resolventwc}. Condition (2)(a) follows directly from
Proposition \ref{p.generator}. As a resolvent, $R(\lambda , A)$
satisfies $\frac{d^k}{d\lambda^k}R(\lambda , A) = (-1)^kk!R(\lambda
, A)^{k+1}$. Interchanging differentiation and integration in the
formula for the Laplace transform  yields
\begin{equation}\label{eq.resest}
\dual{R(\lambda , A)^kx}{y} =  \frac{1}{(k-1)!}\int_0^{\infty}
t^{k-1} \dual{T(t)x}{y} \, dt \,\, ,
\end{equation}
for all $x\in X$ and $y \in Y$. Now let $p$ be a $\tau$-continuous seminorm.
Since $\{ e^{-\alpha t}T(t) \, : \, t \geq 0 \}$ is
$\tau$-equicontinuous, we find a $\tau$-continuous seminorm $q$ 
such that $p(e^{-\alpha t}T(t)x) \leq q(x)$ for all $t\geq 0$
and $x \in X$. Since $\tau$ is consistent, it follows that $\tau = \tau_{\mathfrak{M}}$
for a suitable separating collection of bounded subsets of $Y$. We may thus
assume that $p=p_S$ and $q=p_R$ for certain $S,R \in \mathfrak{M}$.
For $y \in S$, $k \in \CN$ and $\lambda > \alpha$ we obtain from \eqref{eq.resest}
\begin{eqnarray*}
|\dual{(\lambda - \alpha )^kR(\lambda , A)^kx}{y}| & \leq &
\frac{(\lambda - \alpha )^k}{(k-1)!} \int_0^{\infty}
t^{k-1}e^{(\alpha - \lambda )t}|\dual{e^{\alpha t}T(t)x}{y}|\, dt\\
& \leq & \frac{(\lambda - \alpha )^k}{(k-1)!} \int_0^{\infty}
t^{k-1}e^{(\alpha - \lambda )t}
q(x)\, dt\\
& = &  q(x) \,\, .
\end{eqnarray*}
Taking the supremum over $y \in S$, (2)(b) follows.\\
(2) $\Rightarrow$ (1) 
Let $B=A-\alpha$. Since $\tau$ is sequentially complete, if follows from
(2)(b) and the theorem in Section IX.7 of \cite{yosida} that $B$ generates a $\tau$-equicontinuous,
$\tau$-continuous semigroup $\mathbf{T}$ on $X$.
Since $\tau$ is a consistent topology, $L(X, \tau ) \subset L(X, \sigma )$ and
hence $\mathbf{S}$ is a semigroup on $(X,Y)$ (note that $\mathbf{S}$ is
exponentially bounded by Remark \ref{rem1}(1)). Furthermore,
\[ R(\lambda , B) = R\minus\int_0^\infty e^{-\lambda t}S(t)x\, dt \,\, ,\]
where $R\minus\int_0^\infty$ denotes the improper Riemannian
integral with respect to $\tau$. However, since the
map $x \mapsto \dual{x}{y}$ is $\tau$-continuous for
every $y \in Y$, it follows that $\dual{R(\lambda, B)x}{y} =
\int_0^\infty e^{-\lambda t} \dual{S(t)x}{y}\, dt $ for all $x \in
X$ and $y \in Y$. Thus, $\mathbf{S}$ is an integrable semigroup on $(X,Y)$ with generator $B$. 
Now put $T(t) = e^{\alpha t}S(t)$. It is routine to check that $\mathbf{T}$ is
an integrable semigroup with generator $A$. It remains to prove that $\mathbf{T}$
is of type $(M,\omega )$. To that end, consider the rescaled semigroup $\mathbf{T}_{\omega}$.
Note that the generator of $\mathbf{T}_\omega$ is $A-\omega =: C$. 
Now fix $x \in X$ and $y \in Y$. The function $\varphi_{x,y} : t \mapsto \dual{e^{-\omega t}T(t)x}{y}$
has Laplace transform $\dual{R(\lambda , C)x}{y}$. Since $\varphi_{x,y}$ is continuous,
every point $t \geq 0$ is a Lebesgue point of $\varphi_{x,y}$ and we infer from the Post-Widder
inversion formula \cite[Theorem 1.7.7]{abhn} and the equation $\frac{d^k}{d\lambda^k}R(\lambda , C) = (-1)^kk!R(\lambda, C)^{k+1}$ that
\[ \varphi_{x,y}(t) = \lim_{n\to \infty}\dual{[\frac{n}{t}R(\frac{n}{t},C)]^nx}{y} \quad
\forall\, t \geq 0\,\, . \]
However, since $\|\lambda^n R(\lambda , C)^n\| \leq M$, it follows that
$|\dual{e^{-\omega t}T(t)x}{y}| \leq M\norm{x}\cdot \norm{y}$.
Since $Y$ is norming for $X$, it follows that $\mathbf{T}$ has type $(M,\omega )$.
\end{proof}

\begin{rem}
\begin{enumerate}
\item The assumption that $\tau$ is sequentially
complete and consistent is not needed in the implication (1)
$\Rightarrow$ (2) in Theorem \ref{t.hilleyosida}. In the implication
(2) $\Rightarrow$ (1), the sequential completeness is needed to
apply the Theorem from \cite{yosida}, whereas the consistency was used
to ensure that $L(X,\tau) \subset L(X, \sigma )$.
\item The proof of Theorem 3.4 shows that if $\tau$ is sequentially complete and
consistent, then a $\tau$-continuous, quasi-$\tau$-equicontinuous semigroup is
integrable.
\end{enumerate}
\end{rem}

The question remains whether quasi-equicontinuity is a mere technical
assumption in order to prove a Hille-Yosida type theorem or whether
there are interesting cases in which continuity in a certain
topology already implies quasi-equicontinuity. In \cite{komura} it is
proved that on a barreled locally convex space $(X, \tau )$ (recall
that $(X,\tau )$ is called {\it barreled} if every absorbing,
absolutely convex and closed subset of $X$ is a $\tau$-neighborhood
of zero) every $\tau$-continuous semigroup is locally
$\tau$-equicontinuous. The following proposition shows that {\it consistent}
topologies are never barreled, except when the norm topology is
consistent. The special case $(X,Y) = (C_b(E), \mathcal{M}_0(E))$, 
was considered in \cite[Theorem 4.8]{sentilles}.

\begin{prop}
 Let $(X,Y)$ be a norming dual pair and $\tau$ be a consistent topology on $X$. The following
are equivalent.
\begin{enumerate}
 \item $(X, \tau )$ is barreled;
 \item $\tau = \norm{}$ and thus $Y = X^*$.
\end{enumerate}
\end{prop}

\begin{proof}
(1) $\Rightarrow$ (2) If $(X,\tau )$ is barreled, then every weakly
bounded subset of $Y = (X, \tau )'$ is relatively $\sigma'$-compact
and $\tau = \mu $, see \cite[21.4 (4)]{koethe}. However, if every
weakly bounded subset of $Y$ is relatively $\sigma'$-compact, then
$\norm{} = \sup_{\{y: \|y\| \leq 1\}}|\dual{\cdot}{y}|$ is a $\mu$-continuous seminorm
and hence $\mu = \norm{}$. $\quad$ (2)
$\Rightarrow$ (1) Is clear, since every normed space is barreled.
\end{proof}

However also in our general setting there are interesting examples
in which continuity with respect to $\tau_{\mathfrak{M}}$ of a
semigroup on $(X,Y)$ implies quasi-$\tau_{\mathfrak{M}}$-equicontinuity. 
We begin with the following

\begin{lem}\label{l.compact}
Let $(X,Y)$ be a norming dual pair and $\mathfrak{M}$ be a separating collection
of $\sigma'$-compact subsets of $Y$.
Further, let $\Omega$ be a compact Hausdorff space and $F : \Omega
\to L(X,\sigma )$ be strongly $\tau_{\mathfrak{M}}$-continuous. Then
for every $\mathcal{K} \in \mathfrak{M}$ the set
\[ \mathcal{L}_{\mathcal{K}} := \{ F(t)'y \,: \,\, t \in \Omega\, , \, y
\in \mathcal{K}\} \] is $\sigma'$-compact. If for every such
$\mathcal{L}_{\mathcal{K}}$ there exists a set
$\mathcal{K}_0 \in \mathfrak{M}$ such that $\mathcal{L}_{\mathcal{K}} \subset \mathcal{K}_0$, then $\{ F(t)\,:\, t \in \Omega
\}$ is $\tau_{\mathfrak{M}}$-equicontinuous.
\end{lem}

\begin{proof}
We fix $\mathcal{K} \in \mathfrak{M}$ and write for simplicity
$\mathcal{L}$ instead of $\mathcal{L}_{\mathcal{K}}$. Let a
net $z_\alpha = F(t_\alpha )'y_\alpha \in \mathcal{L}$ be given.
Since $\Omega$ is compact, there exists a subnet $t_\beta$ and some
$t_0$ such that $t_\beta \to t_0$ in $\Omega$. Since $\mathcal{K}$
is compact, there is a subnet $y_\gamma$ of $y_\beta$ and an element
$y_0 \in \mathcal{K}$ such that $y_\gamma \weak' y_0$. We claim that
$z_\gamma = F(t_\gamma)y_\gamma \weak' y_0 := F(t_0)y_0$. To see this, let $x
\in X$ be given. Then
\begin{eqnarray*}
|\dual{x}{z_\gamma - z_0}| & \leq & |\dual{F(t_\gamma )x - F(t_0)x}{y_\gamma}| + 
|\dual{F(t_0)x}{y_\gamma - y_0}|\\
& \leq & p_{\mathcal{K}}(F(t_\gamma)x - F(t_0)x) + |\dual{F(t_0)x}{y_\gamma- y_0}| \to 0
\end{eqnarray*}
as $\gamma \to \infty$, by the $\tau_{\mathfrak{M}}$-continuity of
$F(\cdot )x$ and the weak convergence of $y_\gamma$. This shows that
$\mathcal{L}$ is $\sigma'$-compact. We now prove the addendum. If
$\mathcal{L}\subset \mathcal{K}_0 \in \mathfrak{M}$, then
\[ p_{\mathcal{K}}(F(t)x) = \sup_{y\in \mathcal{K}}|\dual{x}{F(t)'y}|
\leq \sup_{y \in \mathcal{L}}|\dual{x}{y}| \leq p_{\mathcal{K}_0}(x) \,\,
\]
for all $x \in X$ and $t \in \Omega$. Hence, if for every
$\mathcal{K} \in \mathfrak{M}$ we find a set $\mathcal{K}_0 \in
\mathfrak{M}$ such that the above holds, it follows that $F(\Omega
)$ is $\tau_{\mathfrak{M}}$-equicontinuous.
\end{proof}

We immediately obtain:

\begin{thm}\label{t.mackey}
Let $(X,Y)$ be a norming dual pair and let $\tau_c:=
\tau_{\mathfrak{C}}$, where $\mathfrak{C}$ is the collection of all
$\sigma'$-compact subsets of $Y$. If $\mathbf{T}$ is a semigroup of 
type $(M,\omega )$ which is $\tau_c$-continuous, then
$\mathbf{T}$ is ${}_{\alpha}$quasi-$\tau_c$-equicontinuous for every
$\alpha > \omega$.
\end{thm}

\begin{proof}
For $\alpha > \omega$ we have $e^{-\alpha t}T(t)x \to 0$ in norm and
hence with respect to the coarser topology $\tau_c$. It follows that
the map
\[
 F : [0,\infty ] \to L(X, \sigma ) \quad,\quad F(t) =
 \left\{\begin{array}{rcl}
            e^{-\alpha t}T(t) &,& t \in [0,\infty )\\
             0 & , & t = \infty
                        \end{array} \right.
\]
is $\tau_c$-continuous. Now the statement follows from Lemma
\ref{l.compact}.
\end{proof}

We note that the topology $\tau_c$ is in general not consistent.
However, it can happen that the Mackey topology $\mu$ coincides with
this topology (\cite{sentilles,conway67} then say $\mu$ is the {\it strong} Mackey 
topology of the pair $(X,Y)$). This is the case if and only if for every
$\sigma'$-compact subset $\mathcal{K}$ of $Y$ also its
$\sigma'$-closed, absolutely convex hull
$\overline{\mathrm{aco}}\,\mathcal{K}$ is $\sigma'$-compact. By
Krein's theorem \cite[24.5 (4)]{koethe}, if $\mathcal{K}$ is
$\sigma'$-compact then $\overline{\mathrm{aco}}\,\mathcal{K}$ is
$\sigma'$-compact if and only if $\overline{\mathrm{aco}}\,\mathcal{K}$ is $\mu'$-complete. In
particular, if $\mu'$ is quasicomplete, i.e. $\mu'$ is complete on
every bounded, $\mu'$-closed subset of $X$, then every
$\mu$-continuous semigroup on $X$ is quasi-$\mu$-equicontinuous.

\begin{cor}\label{c.adjoint}
If $(X,Y)$ is a norming dual pair such that $\mu'$ is quasicomplete,
then every $\mu$-continuous semigroup $\mathbf{T}$ on $(X,Y)$ is
quasi-$\mu$-equicontinuous. In particular
\begin{enumerate}
\item If $\mathbf{T}$ is a norm continuous semigroup on a Banach space $X$,
then its adjoint semigroup $\mathbf{T}^*$ on $X^*$ is $\mu
(X^*,X)$-continuous if and only if it is quasi-$\mu
(X^*,X)$-equicontinuous.
\item If $E$ is a completely regular Hausdorff space such that
$(C_b(E), \beta_0)$ is complete, then every $\mu
(\mathcal{M}_0(E), C_b(E))$-continuous integrable semigroup on
$(\mathcal{M}_0(E), C_b(E))$ is quasi-$\mu (\mathcal{M}_0(E),
C_b(E))$-equicontinuous.
\end{enumerate}
\end{cor}

\begin{proof}
The proof of the general statement was explained above. For (1) we note that
$\mu' = \mu (X, X^*) = \norm{}$ is complete whence every $\mu
(X^*,X)$-continuous adjoint semigroup is quasi-$\mu (X^*,
X)$-equicontinuous. The converse follows from Proposition
\ref{p.taucont} since adjoint semigroups have a $\sigma
(X^*,X)$-densely defined generator and hence, by the Hahn-Banach
theorem, a $\mu (X^*,X)$-densely defined generator. For (2) observe that, 
as a consequence of Grothendieck's completeness theorem \cite[21.9 (4)]{koethe}, the
Mackey topology $\mu (C_b(E), \mathcal{M}_0(E))$ is complete, since there exists a 
complete, consistent topology, namely $\beta_0$, on $C_b(E)$.
\end{proof}

We will now apply Lemma \ref{l.compact} in the context of positive semigroups.
We introduce the following notation. An {\it ordered norming dual pair} 
is a norming dual pair $(X,Y)$ where $X$ is an ordered Banach space with $\sigma$-closed
positive cone $X^+$. Note that in this case the dual cone $Y^+ := \{ y \in
Y\,:\, \dual{x}{y}\geq 0 \,\,\, \forall\, x \in X^+\}$ is $\sigma'$-closed. 
As usual, we call $T\in L(X, \sigma )$ 
{\it positive} if $TX^+ \subset X^+$. Note that in this case also $T'Y^+
\subset Y^+$.

\begin{thm}\label{t.order}
Let $(X,Y)$ be an ordered norming dual pair and $\tau_+$ be the topology of uniform convergence
on the $\sigma'$-compact subsets of $Y^+$. If $\mathbf{T}$ is a positive, $\tau_+$-continuous
semigroup of type $(M,\omega )$ on $(X,Y)$, then $\mathbf{T}$ is ${}_\alpha$quasi- $\tau_+$-equicontinuous for every $\alpha > \omega$.
\end{thm}

\begin{proof}
This follows from Lemma \ref{l.compact} as in the proof of Theorem \ref{t.mackey}, noting that 
for $\alpha > \omega$ and $\mathcal{K} \subset Y^+$ the set $\{e^{-\alpha t}T(t)'y \, : \,
t \geq 0\, , \, y \in \mathcal{K}\}$ is not only compact but also a subset of $Y^+$ by the 
positivity of the operators $T(t)$.
\end{proof}

\section{Applications to Transition semigroups}

In this section, we apply the results of the previous sections to semigroups
on the norming dual pair $(C_b(E), \mathcal{M}_0(E))$. Here, and throughout this section, $E$
will be a completely regular Hausdorff space. The consistent topology
we are interested in is the strict topology $\beta_0$. In order to apply our results,
we need additional information about $\beta_0$ and the dual pair $(C_b(E), \mathcal{M}_0(E))$.
It is well known, see \cite[21.3 (2)]{koethe}, that if $(X,\tau )$ is a locally convex
space then $\tau$ is the topology of uniform convergence on the $\tau$-equicontinuous
subsets of $(X, \tau )'$.  For the strict topology, we have the following
description of the $\beta_0$-equicontinuous subsets of $\mathcal{M}_0(E)$.

\begin{thm}\label{t.sentequi} {\rm (}{\sc Sentilles} \cite[Theorem 5.1]{sentilles}{\rm )}\\
A set  $\mathcal{H} \subset \mathcal{M}_0(E)$ is $\beta_0$-equicontinuous if and only if
(1) $\sup_{\mu \in \mathcal{H}}|\mu |(E) < \infty$ and (2) for every $\eps >0$ there exists a
compact set $K\subset E$ such that
$|\mu |(E\setminus K) \leq \eps$ for all $\mu \in \mathcal{H}$.
\end{thm}

Condition (2) means that $\mathcal{H}$ is a {\it
tight} set of measures. From Theorem \ref{t.sentequi} we infer the following
description of $\beta_0$-equicontinuous sets of linear operators.

\begin{prop}\label{p.equistrict}
Let $\mathcal{S} = \{T_{\alpha} \, : \, \alpha \in I \} \subset
L(C_b(E) , \sigma )$ be a bounded family of operators on $C_b(E)$ with
associated kernels $p_{\alpha}$. The following are equivalent.
\begin{enumerate}
\item $\mathcal{S}$ is $\beta_0$-equicontinuous;
\item given a compact subset $K\subset E$ and $\eps > 0$, there exists a compact
subset $L$ of $E$ such that
\[ |p_{\alpha}|(x, E\setminus L) \leq \eps \quad \forall \, x \in K\,, \, \alpha \in I\,\, .\]
\end{enumerate}
\end{prop}

\begin{proof}
(1) $\Rightarrow$ (2): Let $K \subset E$ be compact. Then
$\mathcal{K} := \{\delta_x \, : \, x \in K\}$ is
$\beta_0$-equicontinuous by Theorem \ref{t.sentequi}. In particular,
$p_{\mathcal{K}}(\cdot )$ is a $\beta_0$-continuous seminorm. Since
$\mathcal{S}$ is $\beta_0$-equicontinuous, we find
a $\beta_0$-continuous seminorm $q$ such that $p_{\mathcal{K}}(T_\alpha f) \leq
q(f)$ for all $f \in C_b(E)$ and $\alpha \in I$. Note that
$q = p_{\mathcal{L}}$ for some $\beta_0$-equicontinuous set
$\mathcal{L}$. We may assume that $\mathcal{L}$ is
$\sigma'$-closed and absolutely convex. But then it follows that
$T_{\alpha}'\delta_x \in \mathcal{L}$ for all $x \in K$. Indeed, if
$T_{\alpha_0}'\delta_{x_0} \not\in \mathcal{L}$ for some $\alpha_0
\in I$ and $x_0 \in K$ then, as a consequence of the Hahn-Banach theorem, we can
strictly separate $T_{\alpha_0}\delta_{x_0}$ from $\mathcal{L}$, i.e. we find
$f \in C_b(E) = (\mathcal{M}_0(E), \sigma')'$ and $\eps >0$ such that
$|\dual{f}{\mu}| + \eps \leq |\dual{f}{T_{\alpha_0}'\delta_{x_0}}|$ for all $\mu \in \mathcal{L}$.
But then  $p_{\mathcal{K}}(f) \geq p_{\mathcal{L}}(f) + \eps$,
a contradiction to the choice of $\mathcal{L}$. Hence, the set $\{
p_{\alpha}(x, \cdot ) \, : \, \alpha \in I\, , \,x\in K \}$ is a
subset of $\mathcal{L}$ and thus $\beta_0$-equicontinuous.
Theorem \ref{t.sentequi} yields (2).\\
(2) $\Rightarrow$ (1): Let $\mathcal{H}$ be a
$\beta_0$-equicontinuous subset of $\mathcal{M}_0(E)$. Then there
exists $C>0$ such that $\norm{\mu} = |\mu |(E) \leq C$ for all
$\mu\in\mathcal{H}$. If we choose $M>0$ such that
$\norm{T_{\alpha}}\leq M$ for all $\alpha \in I$, then
\[ |\dual{f}{T_{\alpha}'\mu}|\leq M\cdot C\cdot \norm{f}
\quad\forall\, f \in C_b(E)\, , \, \alpha \in I\,,\, \mu\in
\mathcal{H}\,\, .
\]
Taking the supremum over $f$ with $\norm{f}_{\infty}\leq 1$, it
follows that $|T_\alpha'\mu|(E) \leq C\cdot M < \infty$ for all
$\alpha \in I$ and $\mu \in \mathcal{H}$. Furthermore, given $\eps
>0$, we find a compact set $K$ such that $|\mu |(E\setminus K ) \leq
\frac{\eps}{2 M}$. By (2), we find $L\subset E$ compact such that
$|p_{\alpha}| (x, E\setminus L) \leq \frac{\eps}{2C}$ for all
$\alpha \in I$ and $x \in K$. It follows that for $\mu \in \mathcal{K}$ and
$\alpha \in I$ we have
\begin{eqnarray*}
|T_{\alpha}'\mu |(E\setminus L) & \leq & \int_K |p_{\alpha}|(x,E\setminus L) \, d|\mu| (x)
+ \int_{K^c} |p_{\alpha}|(x , E\setminus L)| \, d|\mu| (x)\\
& \leq & \frac{\eps}{2C} \norm{\mu} + \frac{\eps}{2 M} M = \eps \,\, .
\end{eqnarray*}
Hence $\mathcal{L} := \{T_{\alpha}'\mu \, : \, \alpha \in I\, , \,
\mu \in \mathcal{H} \, \}$ is $\beta_0$-equicontinuous and thus
$p_{\mathcal{L}}$ is a $\beta_0$-continuous seminorm. However,
$p_{\mathcal{H}}(T_\alpha f) \leq p_{\mathcal{L}}(f)$ for all $f \in
C_b(E)$. Since $\mathcal{H}$ was arbitrary, it follows that
$\mathcal{S}$ is $\beta_0$-equicontinuous.
\end{proof}

Let us recall the following
definition from \cite{sentilles}. A completely regular space $E$ is
called a {\it T-space} if every $\sigma'$-compact set of positive Radon
measures is tight. The celebrated Prokhorov theorem, see \cite{pro},
asserts that every Polish space is a T-space. More generally, every
complete metric space and every locally compact space is a T-space,
see \cite[Theorem 5.4]{sentilles}. If $E$ is an infinite dimensional 
separable Hilbert space endowed with the weak topology, then $E$ is {\it not} a
T-space, cf. \cite[Example I.6.4]{fernique}.

\begin{thm}\label{t.sentilles}
Let $\tau_+$ denote the topology of uniform convergence on the $\sigma'$-compact
subsets of $\mathcal{M}_0^+(E)$.
\begin{enumerate}
 \item $\beta_0 = \tau_+$ iff $E$ is a $T$-space.
 \item If $E$ is a $T$-space and and every measure on $E$ is a Radon measure,
 then $\beta_0 = \mu (C_b(E), \mathcal{M}_0(E))$. In this case, every
$\sigma'$-compact subset of $\mathcal{M}_0(E)$ is tight, i.e. $\beta_0$ 
is the topology of uniform convergence on the $\sigma'$-compact
subsets of $\mathcal{M}_0(E)$.
\end{enumerate}
\end{thm}

\begin{proof}
(1) is \cite[Theorem 5.3]{sentilles}, (2) follows from Theorems 5.8 and 4.5 of that paper.
\end{proof}

We note that {\sc Conway} \cite{conway67} has proved that if $E= [0, \omega_1)$, where $\omega_1$ is
the first uncountable ordinal and $E$ is endowed with the order topology, then $\beta_0$ is not the
Mackey topology of the pair $(C_b(E), \mathcal{M}_0(E))$. However, also in this case
$\beta_0 = \tau_+$ since $E$, being locally compact, is a T-space.

We now come to the main result of this section.

\begin{thm}\label{t.cbmucont}
Let $E$ be a $T$-space and let $\mathbf{T}$ be an integrable semigroup on $(C_b(E), \mathcal{M}_0(E))$.
If there exists a measure $\mu$ on $\mathcal{B}(E)$ which is not a Radon measure, then additionally 
assume that $\mathbf{T}$ is positive.
We denote the kernel associated to $T(t)$ by $p_t$. Consider the following statements
\begin{enumerate}
\item For every $f \in C_b$, the map $(t,x) \mapsto T(t)f(x)$ is continuous.
\item For every $f \in C_b$ we have $T(t)f \to T(s)f$ for $t \to s$ uniformly on the compact
subsets of $E$;
\item $\mathbf{T}$ is a $\beta_0$-continuous semigroup;
\item $\mathbf{T}$ is a quasi-$\beta_0$-equicontinuous, $\beta_0$-continuous semigroup.
\item $\mathbf{T}$ has a $\sigma$-densely defined generator and, given a compact
subset $K\subset E$ and  constants $t_0, \eps > 0$ there exists a compact subset
$L\subset E$ such that
\[ |p_t| (x, E\setminus L) \leq \eps \quad\quad\forall\, x\in K \, , \, t \in [0,t_0] \,\, .\]
\end{enumerate}
Then {\rm (1)} $\Rightarrow$ {\rm (2)} $\Leftrightarrow$ {\rm (3)} $\Leftrightarrow$ {\rm (4)}
$\Leftrightarrow$ {\rm (5)}.
If $E$ is either locally compact or a metric space, all five statements are equivalent.
\end{thm}

\begin{proof}
(1) $\Rightarrow$ (2): Fix $f \in C_b$ and $s\geq 0$. By assumption, for every $\eps > 0$
and $x \in E$ we find $\delta = \delta (s,x)$ and a neighborhood $U = U(s,x)$ of $x$  such that
\[ |T(s)f(x) - T(t)f(y)| < \eps\quad\quad \forall \, (t,y) \in B(s,\delta )\times U \,\, .\]
Now let $K \subset E$ compact. Then $\{s\}\times K \subset
\bigcup_{x \in K} B(s, \delta (s, x)) \times U(s, x)$. As $\{s\}\times K$ is compact in
$[0,\infty )\times E$ we find finitely many $x_1, \ldots , x_n$ and $\delta_i :=
\delta (s, x_i)$ such that $\{s\}\times K \subset \bigcup_{i=1}^n B(s, \delta_i)
\times U(x_i)$. Put $\delta = \min\{\delta_1, \ldots , \delta_n\}$. For $x \in K$,
there exists $i_0$ such that $x \in U(x_{i_0})$. For $|t-s| < \delta$ we have
\[ |T(t)f(x) - T(s)f(x)| \leq |T(t)f(x) - T(s)f(x_i)| + |T(s)f(x_i) - T(s)f(x)|
< 2\eps \, , \]
since $(t,x), (s,x) \in B(s, \delta_{i_0})\times U(x_{i_0})$.
As $x \in K$ was arbitrary, we have $\sup_{x\in K}|T(t)f(x) - T(s)f(x)| \leq \eps$ for
$|t-s|< \delta$. This proves (2).\\
(2) $\Rightarrow$ (1) If $E$ is a metric space, then $(t,x) \mapsto T(t)f(x)$ is
continuous iff it is sequentially continuous. So let $(t_n, x_n) \to (s, x_0)$.
By (2), $T(t_n)f \to T(s)f$ uniformly on the compact set $K=\{x_n \, : \, n\in\CN_0 \}$.
But then $T(t_n)f(x_n) \to T(s)f(x_0)$ follows using the continuity of $T(s)f$.\\
Now assume that $E$ is locally compact. Fix $(s,x_0) \in [0,\infty ) \times
E$ and $f \in C_b$. Since $T(s)f$ is continuous, given $\eps > 0$, there is a neighborhood
$U(x_0)$ such that $|T(s)f(x) - T(s)f(x_0)| < \eps$ for all $x \in U(x_0)$.
It is no restriction to assume that $U(x_0)$ is relatively compact. By (2) we find $\delta > 0$ such that
$|T(t)f(x) - T(s)f(x)| < \eps$ for all $x \in \overline{U(x_0)}$ and all $|t-s|<\delta$.
Thus $|T(t)f(x) - T(s)f(x_0)| < 2\eps$ for all $(t,x) \in B(s, \delta )
\times U(x_0)$. This proves (1).\\
(2) $\Leftrightarrow$ (3): Is clear since $\mathbf{T}$ is locally bounded and
since the strict topology coincides with the compact-open topology on norm-bounded sets.\\
(3) $\Leftrightarrow$ (4): If every measure on $E$ is a Radon measure, then, by Theorem
\ref{t.sentilles} (2), $\beta_0$ is the topology of uniform convergence on the $\sigma'$-compact
subsets of $\mathcal{M}_0(E)$. Hence (3) $\Rightarrow$ (4) follows from Theorem \ref{t.mackey}.
In the other case, $\beta_0$ is the topology of uniform convergence on
the $\sigma'$-compact subsets of $\mathcal{M}_0(E)^+$ by Theorem \ref{t.sentilles} (1). Thus
$\mathbf{T}$ is quasi-$\beta_0$-equicontinuity as a consequence of the positivity of $\mathbf{T}$ 
and Theorem \ref{t.order}. This shows (3) $\Rightarrow$ (4), the converse  implication is trivial.\\
(4) $\Rightarrow$ (5): Follows from Theorem \ref{t.resolventwc} and Proposition
\ref{p.equistrict}.\\
(5) $\Rightarrow$ (4): Is a consequence of Propositions \ref{p.equistrict} and
\ref{p.taucont}.
\end{proof}

\begin{rem}
The assumption in Theorem \ref{t.cbmucont} that $\mathbf{T}$ is an {\it integrable} semigroup 
is only needed for the equivalence (4) $\Leftrightarrow$ (5).
\end{rem}

Theorem \ref{t.cbmucont} can be used to establish that a given
transition semigroup on $C_b(E)$ is $\beta_0$-continuous. In
\cite{belo}, transition semigroups are constructed from the
solutions of parabolic partial differential equations. Here $E$ is a subset of $\CR^d$. 
For such transition semigroups, condition (1) can
easily be verified, as the PDE techniques yield solutions of the PDE which are
continuous in both time and space variables. 
If one follows \cite{cerrai, kuehne, farkas} and prefers to think about semigroups on $C_b(E)$ which have $\tau_{\mathrm{co}}$-continuous orbits, then of course
Condition (2) is satisfied. In the next section, we will show that if $\mathbf{T}$ is the
transition semigroup of a Markov process obtained from solving a stochastic differential
equation, then Condition (5) can often be verified.\\
We note that the crucial assumption in Theorem \ref{t.cbmucont} is
that $\mathbf{T}$ consists of $\sigma (C_b(E), \mathcal{M}_0(E))$-continuous operators.
Under suitable assumption on $E$, e.g. if $E$ is a Polish space,
an operator $T$ on $C_b(E)$ is $\sigma (C_b(E), \mathcal{M}_0(E))$-continuous
if and only if it is a kernel
operator, see \cite{k09a}. If one follows the approach of \cite{belo},
then it is a consequence of the PDE techniques that the operators of
the semigroup are represented by a {\it Green function} and thus are kernel operators.\\
If $E$ is a Polish space and $\mathbf{T}$ is a $\tau_{\mathrm{co}}$-bi-continuous 
semigroup on $C_b(E)$, then it follows from the definition 
of a bi-continuous semigroup that every operator $T(t)$ is sequentially 
$\tau_{\mathrm{co}}$-continuous on normbounded sets and hence sequentially $\beta_0$-continuous.
By \cite[Corollary 8.4]{sentilles} it follows that $T(t) \in L(C_b(E), \beta_0) \subset 
L(C_b(E), \sigma )$. {\sc Farkas} \cite[Example 3.9]{farkas09} has given an example of a
$\tau_{\mathrm{co}}$-bi-continuous semigroup which does not consist of $\beta_0$-continuous operators
and is thus, in particular, not locally $\beta_0$-equicontinuous.
In that example, $E = [0, \omega_1)$ with the order topology, where $\omega_1$ is the first
uncountable ordinal. Note that since $[0, \omega_1)$ is locally compact
and hence a T-space, it follows from Theorem \ref{t.cbmucont} that every positive 
$\beta_0$-continuous semigroup of operators in $L(C_b(E), \sigma )$ is 
quasi-$\beta_0$-equicontinuous.

\section{Examples}

\subsection{The Case $E = \CN$} If $E = \CN$ endowed with the
discrete topology, then $C_b(E) = \ell^\infty$ and 
$\mathcal{M}_0(E)= \ell^1$. Thus in this case, $\mathcal{M}_0(E)$ is
the predual of $C_b(E)$. The weak topology $\sigma = \sigma (C_b(E),
\mathcal{M}_0(E))$ is merely the $\mathrm{weak}^*$-topology of
$\ell^\infty$ whereas the weak topology $\sigma' = \sigma
(\mathcal{M}_0(E), C_b(E))$ is the weak topology (in the Banach
space sense) of $\ell^1$. A bounded operator $T$ on $\ell^\infty$ is
$\sigma$-continuous if and only if it is the adjoint of a bounded
operator on $\ell^1$. We now have the following result.

\begin{prop}
If $E= \CN$ endowed with the discrete topology, then every semigroup
$\mathbf{T}$ on $(C_b(E),\mathcal{M}_0(E))$ which is
$\sigma$-continuous at 0 is $\beta_0$-continuous and quasi-$\beta_0$-equicontinuous.
\end{prop}

\begin{proof}
If $\mathbf{T}$ is a semigroup on $(\ell^{\infty},\ell^1)$ then
$\mathbf{T}=\mathbf{S}^*$ for some semigroup $\mathbf{S}$ on
$\ell^1$. Now $\mathbf{T}$ is $\sigma
(\ell^\infty,\ell^1)$-continuous (at 0) if and only if $\mathbf{S}$ is
$\sigma (\ell^1,\ell^\infty )$-continuous (at 0). It is well known (but
also follows from Theorem \ref{t.resolventwc} and Corollary
\ref{c.weakequalstrong}) that a semigroup of bounded operators on a
Banach space $X$ which is weakly continuous at 0 is already
$\norm{}$-continuous. It follows that $\mathbf{T}$ is
$\sigma(\ell^\infty,\ell^1)$-continuous. In particular, $t \mapsto
T(t)f(n)$ is continuous for every $n \in \CN$ and every $f \in
C_b(\CN )$. However, since every compact subset of $\CN$ is finite,
it follows that $t \mapsto T(t)f$ is $\tau_{\mathrm{co}}$-continuous
and hence $\beta_0$-continuous for every $f \in C_b(\CN )$. Since
$\CN$ is locally compact and since every measure on $\CN$ is a Radon measure, 
the assertion follows from Theorem \ref{t.cbmucont}.
\end{proof}

\subsection{The Sorgenfrey line}

Let us consider the real line $\CR$ endowed with the Sorgenfrey
topology $\tau_s$, i.e. $\tau_s$ is generated by the intervals of
the form $[a,b)$. It follows that the Borel $\sigma$-algebra of
$(\CR, \tau_s )$ is the usual Borel $\sigma$-algebra of $\CR$
(with the metric $|\,\cdot \, |$). It is well known that every
compact subset of $(\CR, \tau_s )$ is countable. However, as the
example $\{1- \frac{1}{n}\, :\, n \in \CN\}$ shows, not every
countable subset of $\CR$ is compact in $(\CR, \tau_s )$. It follows
that every Radon measure on $(\CR, \tau_s )$ is concentrated on a
countable set. Hence $\mathcal{M}_0(\CR,\tau_s ) = \ell^1(\CR )$,
the space of all discrete measures on $\CR$. A function $f$ on $\CR$
is continuous with respect to $\tau_s$ if and only if it is
right-continuous. Thus in this situation $(C_b(\CR, \tau_s ),
\mathcal{M}_0(\CR,\tau_s )) = (C_r(\CR), \ell^1(\CR ))$, where
$C_r(\CR)$ denotes the space of all bounded, right continuous
functions on $\CR$.\\
Consider the shift semigroup $\mathbf{T}$ given by $T(t)f(x) =
f(x+t)$. Then $\mathbf{T}$ is a positive contraction semigroup on 
$(C_r(\CR),\ell^1(\CR ))$. However, $\mathbf{T}$ is not integrable.
Indeed, it is easy to see that the Laplace transform of $\mathbf{T}$ 
is given by
\[ R(\lambda )f(x) = e^{\lambda x}\int_x^\infty e^{-\lambda y}f(y)\, dy \,\, . \]
But this operator is not $\sigma$-continuous since
$R(\lambda )^*\delta_0 = e^{-\lambda t}\one_{(0,\infty )}dt \not\in \ell^1(\CR )$.
Furthermore, $\mathbf{T}$ is $\sigma$-continuous at $0$ but not $\sigma$-continuous. 
Since a continuous function
is uniformly continuous on compact sets, it follows that
$T(t)f \stackrel{\tau_{\mathrm{co}}}{\to} f$ as $t \downarrow 0$ for every $f \in C_r(\CR)$. 
Hence $\mathbf{T}$ is $\beta_0$-continuous at $0$. 

\subsection{Solutions of stochastic differential equations}
If $E$ is a Banach space, then some Markov processes are obtained as
solutions of stochastic differential equations, see e.g. \cite{gk01, gvn03}.
The transition semigroup of such a Markov process is defined as follows. 
If $X(t,x)$ denotes the solution of the stochastic differential equation 
with initial datum $x \in E$, then one defines $T(t)f(x) = \mathbb{E}(f(X(t,x)))$. 
It is natural to ask for a condition for $\mathbf{T}$ to be $\beta_0$-continuous
in terms of properties of the map $(t,x) \mapsto X(t,x)$. In applications, it
frequently happens that $X(t,x) \in L^p(\Omega, E)$, where $\Omega$ is
the underlying probability space and $1\leq p < \infty$, and the map $(t,x)
\mapsto X(t,x)$ is continuous. The following theorem shows that in this case
the semigroup $\mathbf{T}$ is indeed $\beta_0$-continuous. However, this
result remains true in a even more general setting.

If $(\Omega, \mathscr{F}, P)$ is a probability space and $(E, \rho )$ is a 
complete metric space, then $L^0(\Omega , E)$ denotes the space of
all strongly measurable maps $X : \Omega \to E$ (modulo equality
$P$-almost everywhere) endowed with the topology of convergence in measure.

\begin{thm}
Let $(\Omega, \mathscr{F}, P)$ be a probability space, $(E, \rho)$ be a complete metric
space and $X : [0,\infty ) \times E \to L^0(\Omega , E)$ be a continuous map. 
Define $T(t)f(x) = \mathbb{E}(f(X(t,x))$
for $f\in C_b(E)$. Then, for every $t_0>0$, the set  $\{T(t) \, : 0 \leq t \leq t_0\}$ is a $\beta_0$-equicontinuous family
of operators on $(C_b(E), \mathcal{M}_0(E))$. If $(T(t))_{t\geq 0}$
is a semigroup, then it is $\beta_0$-continuous and quasi-$\beta_0$-equicontinuous.
\end{thm}

\begin{proof}
Consider the map $\Phi : L^0(\Omega , E ) \to \mathcal{M}_0(E)$ given by
$\Phi (X ) = \mu_X$, where $\mu_X$ denotes the distribution of $X$. 
Note that for $X \in L^0(\Omega , E)$ the distribution $\mu_X$ is indeed
a Radon measure since $X$ has separable range. The map $\Phi$ is continuous.
Indeed, if $X_n \to X$ in measure then, passing to a subsequence, we have
$X_n \to X$ almost everywhere. But then 
\[ \dual{f}{\mu_{X_n}} = \int_{\Omega} f(X_n) \, dP \to \int_{\Omega} f(X)\, dP =
 \dual{f}{\mu_X} \,\, .
\]
Thus, every subsequence of $\Phi (X_n)$ has a subsequence which
converges to $\Phi (X)$ with respect to $\sigma (\mathcal{M}_0(E),C_b(E))$. 

Since $T(t)f(x) = \dual{f}{\mu_{X(t,x)}}$, the continuity of $x \mapsto
X(t,x)$ for fixed $t$ implies that  every operator $T(t)$ 
maps $C_b(E)$ into $C_b(E)$. It follows from the joint continuity of
$X(\cdot , \cdot )$ that for every $t_0>0$ and every compact set
$K\subset E$ the set $\{ X(t,x) \, : \, 0\leq t \leq t_0 \, , x \in
K \}$ is compact in $L^0(\Omega , E)$. Hence the set $\{\mu_{X(t,x)}\, : \, 0
\leq t \leq t_0\, , \, x \in K \}$ is $\sigma'$-compact in
$\mathcal{M}_0(E)$. Since $E$ is a complete metric space and hence a
T-space, it follows that the latter set is tight. Hence Proposition
\ref{p.equistrict} implies that the set $\{T(t)\, :\, 0 \leq t \leq t_0\} $ is 
$\beta_0$-equicontinuous. In particular, every single operator $T(t)$ is
$\beta_0$-continuous and hence an element of $L(C_b(E), \sigma )$.
Now assume that $(T(t))_{t\geq 0}$ is a semigroup. Since $t \mapsto
X(t,x)$ is continuous, it follows that $t \mapsto T(t)f(x)$ is
continuous for every $x \in E$ and hence $(T(t))_{t\geq 0}$ is
$\sigma$-continuous. Taking Remark \ref{rem1} into account, it follows from
Proposition \ref{t.resolventwc} that $\mathbf{T}$ has a $\sigma$-densely defined generator. 
The claim follows from Theorem \ref{t.cbmucont}.
\end{proof}

\end{document}